\documentstyle[12pt,twoside]{amsart}

\let \cedilla =\c
\renewcommand{\c}[0]{{\mathbb C}}  

\renewcommand{\o}[0]{{\cal O}} 
\newcommand{\z}[0]{{\mathbb Z}}

\renewcommand{\r}[0]{{\mathbb R}} 

\renewcommand{\a}[0]{{\mathbb A}}

\newcommand{\p}[0]{{\mathbb P}}

\newcommand{\qtq}[1]{\quad\mbox{#1}\quad}
\newcommand{\spec}[0]{\operatorname{Spec}}

\newcommand{\rank}[0]{\operatorname{rank}}
\newcommand{\mult}[0]{\operatorname{mult}}

\newcommand{\supp}[0]{\operatorname{Supp}}

\newsymbol\subsetneq 2328
\newsymbol\onto 1310

\newsymbol\twoheadrightarrow 1310

\newtheorem{thm}{Theorem}[section]
\newtheorem{question}[thm]{Question}
\newtheorem{lem}[thm]{Lemma}
\newtheorem{cor}[thm]{Corollary}

\newtheorem{prop}[thm]{Proposition}

\theoremstyle{definition}
\newtheorem{defn}[thm]{Definition}

\newtheorem{say}[thm]{}
\newtheorem{exmp}[thm]{Example}

\newtheorem{rem}[thm]{Remark}          
\newtheorem{ack}{Acknowledgments}        
\newtheorem{notation}[thm]{Notation}

\theoremstyle{remark}

\begin{document}
\bibliographystyle{amsplain}

\title{Effective Nullstellensatz for Arbitrary Ideals}
\author{J\'anos Koll\'ar}

\maketitle
\tableofcontents

\newcommand{\tdeg}[0]{\operatorname{arith-deg}}    
\newcommand{\gdeg}[0]{\operatorname{geom-deg}}    
\newcommand{\dist}[0]{\operatorname{dist}}   
\newcommand{\ann}[0]{\operatorname{Ann}}   
\newcommand{\norm}[1]{||#1||}   
\newcommand{\len}[0]{\operatorname{length}}   
\newsymbol\ncap 1374  
\newsymbol\dcap 1265

\section{Introduction}

Let $X,Y\subset \p^n$ be closed irreducible subvarieties and 
$Z_i$  the irreducible components of $X\cap Y$.  One variant of the
theorem  of B\'ezout (cf. \cite[8.4.6]{fulton}) says that 
$$
\sum_i\deg Z_i\leq \deg X\cdot \deg Y.
$$ This result holds without any restriction on the dimensions of
 $X,Y, Z_i$ and it can be easily generalized to the case when
$X_1,\dots,X_s$ are arbitrary subschemes of $\p^n$ and the $Z_i$
are the reduced irreducible components of $X_1\cap\dots\cap X_s$.

It is frequently of interest to study finer algebraic or metric
properties of intersections of varieties. In recent years
considerable attention was paid to the case when the $X_i$ are all
hypersurfaces, in connection with the effective versions of
Hilbert's  {\it Nullstellensatz}.  Assume that we have  polynomials
$f_1,\dots,f_s\in
\c[x_1,\dots,x_n]$ of degrees $d_i=\deg f_i$. There are three
related questions one can ask about the intersection of the
hypersurfaces $(f_i=0)$, in each case attempting to  minimize a
bound $B(d_1,\dots,d_s)$.

\begin{description}
\item[Algebraic B\'ezout version]\cite{brownaw-ppp}
 Find prime ideals $P_j\supset (f_1,\dots,f_s)$ and natural numbers
$a_j$ such that
$$
\prod_jP_j^{a_j}\subset (f_1,\dots,f_s)
\qtq{and} \sum_j a_j\deg P_j\leq B(d_1,\dots,d_s).
$$

\item[Effective Nullstellensatz version]\cite{null}
  If the $f_i$ have no common zeros in $\c^n$, find polynomials
$g_i$  such that
$$
\sum_i f_ig_i=1
\qtq{and} \deg (f_ig_i)\leq B(d_1,\dots,d_s).
$$

\item[\L ojasiewicz inequality version] \cite{jks}
 Fix a metric on $\c^n$ and let $Z$ be the intersection of the
hypersurfaces $(f_i=0)$.  Prove that if $x$ varies in a bounded
subset  of $\c^n$ then
$$
\dist(Z,x)^{B(d_1,\dots,d_s)}\leq C\cdot \max_i|f_i(x)|\qtq{for some
$C>0$.}
$$
\end{description}

The optimal value of the bound $B(d_1,\dots,d_s)$ is known in
almost all cases. If we assume that $d_i\geq 3$ for every $i$, then 
$$ B(d_1,\dots,d_s)=d_1\cdots d_s,
$$ is best possible for $s\leq n$. (See 
\cite[1.5]{null} for the case
$s>n$.)

The algebraic B\'ezout version  is also called the {\it prime power
product} variant of  the Nullstellensatz. 
\medskip

The aim of this paper is to consider these problems in case the 
$f_j$ are replaced by arbitrary ideals. The first step in this
direction was taken in \cite{sombra}. His methods can deal with 
special cases of the above problems if the ideals are
Cohen--Macaulay.
 \L ojasiewicz--type inequalities for arbitrary analytic sets were
studied in the  works of Cygan, Krasi\'nski and Tworzewski, see
especially
\cite{twor, cygan, ckt}. Although they consider the related problem
of separation exponents, their proof can easily be modified to give
a general \L ojasiewicz inequality for reduced subschemes. 

My proofs grew out of an attempt to understand their work in
algebraic terms. This leads to 
 a general \L ojasiewicz inequality in the optimal form and to an
effective Nullstellensatz with a slightly worse bound. In the
algebraic B\'ezout version my results are considerably weaker. It
should be noted, however, that the straightforward generalization
of the algebraic B\'ezout version fails to hold
(\ref{algbez.false.exmp}).

All three  of these results can be formulated for arbitrary ideals,
but for simplicity here I state them  for unmixed ideals. ($I$ is
called {\it unmixed} if all   primary components of $I$ have the
same dimension.) These are the ideals that correspond to the usual
setting of intersection theory.
   For such ideals the degree of $I$  gives a good generalization
of the degree of a hypersurface. The precise versions for arbitrary
ideals are stated in   (\ref{bez.thm}), (\ref{eff.null.thm}) and 
(\ref{lojas.thm}). 

\begin{thm}[Algebraic B\'ezout theorem]\label{main.intr.thm}
  Let $K$ be any field and \linebreak 
$I_1,\dots,I_m$ unmixed   ideals in $K[x_1,\dots,x_n]$.  Then there
are   prime ideals $P_j\supset (I_1,\dots,I_m)$ and natural numbers
$a_j$ such that
\begin{enumerate}
\item  $\prod_j P_j^{a_j}\subset (I_1,\dots,I_s)$, and
\item $\sum_j a_j\leq n\cdot \prod_i\deg I_i$.
\end{enumerate}
\end{thm}

\begin{thm}[Effective Nullstellensatz]\label{eff.null.intr.thm} Let
$K$ be any field and  \linebreak 
$I_1,\dots,I_m$ unmixed   ideals in
$K[x_1,\dots,x_n]$.  The following are equivalent:
\begin{enumerate}
\item  $I_1,\dots,I_m$ have no common zero in $\bar K^n$.
\item There are polynomials $f_j\in I_j$  such that
$$
\sum_jf_j=1\qtq{and} \deg f_j
\leq (n+1)\cdot 
\prod_i\deg I_i.
$$
\end{enumerate}
\end{thm}

\begin{thm}[\L ojasiewicz inequality]\label{lojas.intr.thm} (cf.
\cite{cygan, ckt}) Let $I_1,\dots,I_m$ be unmixed   ideals in
$\c[x_1,\dots,x_n]$ and $X_1,\dots,X_m\subset \c^n$ the
corresponding  subschemes. Let $f_{ij}$ be generators of  $I_i$. 
Then for every bounded set $B\subset \c^n$ there is a $C>0$ such
that for every $x\in B$,
$$
\dist(X_1\cap \cdots \cap X_m,x)^{\prod_i\deg I_i}\leq C\cdot
\max_{ij}|f_{ij}(x)|.
$$
\end{thm}

The difference between the geometric and algebraic versions of the
B\'ezout theorem can be seen already in the case when an
irreducible variety is intersected with a hyperplane.

\begin{exmp}\label{algbez.false.exmp}
 Pick coordinates $u,v$ in $\c^2$ and $x,y,z,s$ in $\c^4$. For
every odd $n$ consider the morphism
$$ F_n:\c^2\to \c^4\qtq{given by} F_n(u,v)=(u^n,u^2,uv,v).
$$ Let $S_n\subset \c^4$ be the image of $F_n$. It is easy to see
that
$\deg S_n=n$ and 
 the ideal of $S_n$ in $\c[x,y,z,s]$ is 
$$ I_n=(x^2-y^n,z^2-ys^2, z^n-xs^n, xz-ys^{\frac{n+1}{2}}, 
xs-y^{\frac{n-1}{2}}z).
$$ Let us intersect $S_n$ with the hyperplane $(s=0)$ to get a
curve $C_n$. Set theoretically, the intersection is the image of
$f_n:\c\to \c^4$ given by $f_n(u)=(u^n,u^2,0,0)$ and its
 ideal  is
$$ J_n=(x^2-y^n,z,s).
$$ On the other hand,
$$ (I_n,s)=(x^2-y^n,z^2,  xz, y^{\frac{n-1}{2}}z,s),
$$ and we see that, as a vectorspace,
$$ J_n/(I_n,s) \cong \langle z,yz,\dots, y^{\frac{n-3}{2}}z\rangle.
$$ Let $m=(x,y,z,s)$ be the ideal of the origin. There are  two
minimal ways of writing an algebraic B\'ezout form of this example:
$$ J_n^2\subset (I_n,s)\qtq{and} m^{\frac{n-1}{2}}\cdot J_n\subset
(I_n,s).
$$ Taking degrees we get
$2\deg J_n=2n>n=\deg I_n$ and $\frac{n-1}{2}\deg m+\deg J_n=
\frac{n-1}{2}+n>n=\deg I_n$.
\end{exmp}

This example illustrates the nature of the difficulties, but it
does not seem to give   pointers as to the general shape of the
theory. Unfortunately, I do not have any plausible conjectures
about what happens in general. As in \cite{null}, the effect of
embedded primes seems small, but the correct way of estimating it
is still elusive.

Instead, I  approach the question as follows. There are many
different varieties  $S^{\lambda}_n\subset\c^4$ whose intersection
with  the hyperplane
$(s=0)$ is $C_n$.  (For instance,  pick polynomials $f(u,v),g(u,v)$
with  no common zero and let
$S^{f,g}_n$ be the image of $(u,v)\mapsto
(u^n,u^2,vf(u,v),vg(u,v))$.) Each  $S^{\lambda}_n$ gives an ideal
$I^{\lambda}_n$ and one can ask about all the quotients
$$ J_n/(I^{\lambda}_n,s).
$$ It turns out that their length is  bounded independent of
$\lambda$ and it is not too big. The main lemma of \cite[3.4]{null}
is a formalization of this observation using local cohomology groups
in some special cases.

This paper develops another  approach to this problem, going
back to \cite{cayley}. For any space curve $C\subset \p^3$ Cayley
considered  all
 cones defined by $C$ with a variable point $p\in \p^3$ as  
vertex. These cones can be encoded as one equation on the
Grassmannian of lines in $\p^3$.  More generally, for any pure
dimensional subscheme  $Y^d\subset \c^n$ (or  for
any  pure dimensional algebraic cycle on
$\c^n$)  consider the ideal  $I^{ch}(Y)$  generated by all cones
defined by $Y$ with   a variable $(n-d-1)$-dimensional linear space
as its vertex. Following \cite{ds95}, it is called  the {\it ideal
of Chow equations} (\ref{con-id.def}).  It turns out that this ideal
 controls the length of the embedded components of any
intersection. With this observation at hand, the rest of the
arguments turn out to be not very complicated.

Section 2  reviews some basic facts about algebraic cycles and
their intersection theory on $\c^n$. Section 3 collects known
results about integral closures of ideals.

The ideal of Chow equations is defined  and studied in section 4.
The connection between the ideal of Chow equations and intersection
theory is established in section 5.

Finally the main results are proved in sections 6 and 7.
\medskip

Another approach to such theorems is  to reduce them to the
hypersurface case. If $X\subset \c^n$ is a subscheme of degree
$m$ then,  set theoretically, $X$ can be defined by  degree $m$
equations. This gives reasonable bounds for each problem,  roughly
like $(\max_i\{\deg I_i\})^n$. For many ideals of about the same
degree this is close to the optimal bound for the Nullstellensatz,
but  it is considerably worse  in general. For the algebraic
B\'ezout version this method  and \cite{brownaw-ppp} gives a bound
in the original form taking into account the degrees of the $P_j$.

A modified version of this idea is to reduce everything to
intersecting with the diagonal and then using the methods of
\cite{null} directly. This gives  $3^n$-times the optimal bounds.
If, however, the quotients
$K[x_1,\dots,x_n]/I_j$  (or more precisely,  their homogenizations)
are  Cohen-Macaulay,  then the methods of
\cite{sombra} give better bounds. The factor $(n+1)$ in 
(\ref{eff.null.intr.thm}) can be replaced by 2.

\begin{ack}  I   thank  J.\ Johnson and P.\ Roberts for their help
with many of the computations. I am greatful for the comments of
A.\ Beauville, P.\ Philippon  and P. Tworzewski 
 and for several e--mails of B.\ Sturmfels  which  helped to
eliminate many mistakes.

Partial financial support was provided by  the NSF under grant
number  DMS-9622394. 
\end{ack}

\section{Intersection   of cycles on $\a^n$}

\begin{defn}\label{cyc.def} Let $Y$ be a scheme. An {\it algebraic
cycle} on $Y$ is a formal linear combination of reduced and
irreducible subschemes $Z=\sum a_i[Z_i]$,
$a_i\in \z$. (I do not assume that the $Z_i$ have the same
dimension.) The cycles form a free Abelian group $Z_*(Y)$. The
subgroup generated by all reduced and irreducible subschemes of
dimension $d$ is denoted by $Z_d(Y)$.

If $Y$ is proper and $L$ is a line bundle on $Y$ then one can define
the $L$-degree of a cycle 
$$
\deg_LZ:=\sum a_i(Z_i\cdot L^{\dim Z_i}),
$$ where $(Z_i\cdot L^{\dim Z_i})$ denotes the top selfintersection
number of the first Chern class of $L|_{Z_i}$. The function
$Z\mapsto \deg_LZ$ is linear.

Let $Y$ be a scheme with a compactification $Y\subset\bar Y$ and
assume that $L$ is the restriction of a line bundle
$\bar L$ from $\bar Y$ to $Y$.  For a cycle $Z=\sum a_iZ_i$ on $Y$
set $\bar Z=\sum a_i\bar Z_i$ where $\bar Z_i$ is the closure of
$Z_i$ in $\bar Y$.  Then one can define the degree of a cycle
$Z=\sum a_iZ_i$ on $Y$ by
$$
\deg_LZ:=\deg_{\bar L}\bar Z.
$$ It is important to note that this depends on the choice of $\bar
Y$ and $\bar L$. I use this version of the degree
 only for the pair
$Y=\a^n$ and $\bar Y=\p^n$.
\end{defn}

\begin{defn}\label{ncap.def}
 Let $X$ be a scheme and $D$ an effective  Cartier divisior on $X$.
Let $[Z]\in Z_d(X)$ be an irreducible $d$-cycle on $X$. Define
$Z\ncap D \in Z_*(X)$ as follows.
\begin{enumerate}
\item If $Z\subset \supp D$ then set $[Z]\ncap D:=[Z]\in Z_d(X)$.
\item If $Z\not\subset \supp D$ then $D|_{Z}$ makes sense as a
Cartier divisor. Set $[Z]\ncap D:=[D|_{Z}]\in Z_{d-1}(X)$.
\end{enumerate} This definition can be extended to $Z_*(X)$ by
linearity.

If $f$ is a defining equation of $D$ then I also use  
$Z\ncap f$ to denote $Z\ncap D$.

It  should be emphasized that this definition is not at all well
behaved functorially. While it is well defined on cycles, it is not
well defined on the Chow group. Furthermore, if $D_1,D_2$ are two
Cartier divisors then in general
$$ (Z\ncap D_1)\ncap D_2\neq (Z\ncap D_2)\ncap D_1.
$$ (For instance let $X=\a^2, Z=(y-x^2=0), D_1=(x=0)$ and
$D_2=(y=0)$.)
\end{defn}

\begin{lem}\label{ncap-deg.lem} 
 Let $L$ be an ample line bundle on $X$,
  $D$  a section of $L^{\otimes d}$ and $Z$ an effective cycle on
$X$. Then
\begin{enumerate}
\item
$\deg_L(Z\ncap D)\leq d\cdot \deg_L(Z)$.
\item  If $X$ is proper, $d=1$ and all
 the components of $Z$ have positive dimension then $\deg_L(Z\ncap
D)= \deg_L(Z)$.
\end{enumerate}
\end{lem}

Proof. By linearity it is sufficient to check this when $Z$ is an
irreducible and reduced subvariety. If $Z\subset \supp D$ then
$\deg_L(Z\ncap D)= \deg_L(Z)$, and otherwise
$\deg_L(Z\ncap D)\leq d\cdot \deg_L(Z)$  with equality holding if
 $X$ is proper and
 $\dim Z\geq 1$ by the usual B\'ezout theorem.\qed

\begin{say} One would like to define $Z_1\ncap Z_2$ for any two
cycles
$Z_i$ on a scheme $X$. As usual, this is reduced to intersecting
$Z_1\times Z_2$ with the diagonal  $\Delta\subset X\times X$.
Traditional intersection theory  works if $X$ is smooth since in
this case $\Delta\subset X\times X$ is a local complete intersection
(cf. \cite[Chap.8]{fulton}). The usual intersection product 
$Z_1\cdot Z_2$     is then a cycle of the expected dimension
$d=\dim Z_1+\dim Z_2-\dim X$. If $\dim (Z_1\cap Z_2)=d$ then
$Z_1\cdot Z_2$  is well defined as a cycle, but if 
$\dim (Z_1\cap Z_2)>d$ then $Z_1\cdot Z_2$ is defined only as a
rational equivalence class inside $\supp(Z_1\cap Z_2)$.

Here I follow the path of \cite{stu-vo, vogel} and try to define
$Z_1\ncap Z_2$ as a well defined cycle which may have components of
different dimension.  If  $X=\p^n$, the $Z_i$ are pure dimensional
and 
$d\geq 0$ then $Z_1\ncap Z_2$ is a cycle 
 such that
$$
\deg (Z_1\ncap Z_2)=\deg Z_1\cdot\deg Z_2.
$$ The cases when $d<0$ were not considered to have much meaning
traditionally. \cite{twor} realized that  the definition is
meaningful and gives an interesting invariant.

The construction of (\ref{ncap.def}) needs $\Delta$ to be a global
complete  intersection. Unfortunately this happens very rarely. The
only such example that comes to mind is $X=\a^n$, or more
generally, any scheme $X$ which admits an \'etale map to
$\a^n$. For simplicity of exposition, I work with $X=\a^n$.
Homogenity considerations can then be used to define $\ncap$ for a
few  other interesting cases, most importantly for
$X=\p^n$.
\end{say}

\begin{defn}[Vogel--Tworzewski cycles]\label{vt-cyc.def}
 Let $X_i=\sum_j a_{ij}X_{ij}$
 be  effective cycles on $\a^n$ for $i=1,\dots,s$. We would like to
define a cycle which can reasonably be called the intersection of
these  cycles. This is done as follows.

Choose an identification $\a^{ns}=\a^n\times\cdots\times\a^n$.
Using this identification  define
$$
\prod_{i=1}^s X_i:=\sum_{j_1,\dots,j_s}
\left(\prod_{i=1}^sa_{ij_i}\right)
\left(\prod_{i=1}^sX_{ij_i}\right)
$$ as a cycle in $Z_*(\a^{ns})$. 

Let $\Delta\subset \a^n\times\cdots\times\a^n$ denote the 
diagonal. Each coordinate projection
$$
\Pi_r:\a^n\times\cdots\times\a^n\to \a^n\qtq{(onto the $r$th
factor)}
$$ gives an isomorphism $\Pi_r:\Delta\cong \a^n$ which is
independent of $r$. 

Let ${\cal L}:=(L_1,\dots,L_{n(s-1)})$ be  an ordered set of
hyperplanes in $\a^{ns}$ such that their intersection is $\Delta$.
Set
$$ (X_1\ncap\cdots\ncap X_s, {\cal L}):=
\left(\prod_{i=1}^s X_i\right)\ncap L_1\ncap \cdots\ncap L_{n(s-1)},
$$ where the right hand side means that we first intersect with
$L_1$,
 then with $L_2$ and so on. To be precise, the right hand side is
in $Z_*(\a^{ns})$, but every irreducible component of it is
contained in $\Delta$. Thus it
 can be viewed as a cycle in $Z_*(\Delta)$ and so it can be
identified with a cycle in $Z_*(\a^n)$ using any of the projections
$\Pi_r$.

$(X_1\ncap\cdots\ncap X_s, {\cal L})$ is called an {\it
intersection cycle} of $X_1,\dots,X_n$. Any of these cycles is
denoted by
$X_1\ncap\cdots\ncap X_s$. 

It should be emphasized that 
$X_1\ncap\cdots\ncap X_s$ is not a well defined cycle since it
depends on the choice of ${\cal L}$.  In the papers \cite{vogel,
vangast}  the $L_i$ are chosen generic and then 
$(X_1\ncap\cdots\ncap X_s, {\cal L})$ is well defined as an element
of a suitable Chow group.  We would like to get a cycle which is
defined over our field $K$. As long as $K$ is infinite, a general
choice of the $L_i$ would work but there are some problems when $K$
is finite. (It is for such reasons that \cite{brownaw-ppp} does not
work for all finite fields.)
 Furthermore, in our applications it is sometimes advantageous to
make a  special choice of the $L_i$. For these reasons I allow any
choice of the $L_i$.  The price we pay is that 
 even the degree of
$(X_1\ncap\cdots\ncap X_s, {\cal L})$ depends on the $L_i$. This,
however, does not seem to cause  problems in the applications.
\end{defn}

We obtain the following B\'ezout type inequality.

\begin{thm}\label{bezout.thm}
 Let $X_1,\dots,X_s$ be effective  cycles on $\a^n$. Then 
$$
\deg (X_1\ncap\cdots\ncap X_s)\leq \prod_j\deg X_j.
$$
\end{thm}

Proof. $\deg \prod_{i=1}^s X_i=\prod_{i=1}^s \deg X_i$ and cutting
with a hyperplane does not increase the degree by
(\ref{ncap-deg.lem}).\qed

\begin{defn}[Refined intersection cyle] \label{vog.cyc.def} 

Let $K$ be an infinite field.
For a scheme $Y$ let $B(Y)$ denote all subvarieties of $Y$ which can
be obtained by repeatedly taking irreducible components and their
intersections. For general $\cal L$ we can write
$$ (X_1\ncap\cdots\ncap X_s, {\cal L})=\sum a_i[Z_i({\cal L})],
$$ where the $Z_i({\cal L})$ depend algebraically on $\cal L$. For
each $Z_i({\cal L})$ there is a smallest 
$W\in B(X_1\cap\cdots\cap X_s)$ such that $Z_i({\cal L})\subset W$
for every general choice of ${\cal L}$. For each 
 $W\in B(X_1\cap\cdots\cap X_s)$, the sum of these cycles gives a
well defined element of the Chow group $A_*(W)$.  This cycle is
denoted by $(X_1\dcap\cdots\dcap X_s,W)$.  Thus we obtain a refined
intersection cycle
$$ X_1\dcap\cdots\dcap X_s := \sum_{W\in B(X_1\cap\cdots\cap
X_s)}(X_1\dcap\cdots\dcap X_s,W)
$$ If $Z\subset X_1\cap\cdots\cap X_s$ is a connected component then
$$
\deg (X_1\dcap\cdots\dcap X_s,Z):=\sum_{W\subset Z}
\deg (X_1\dcap\cdots\dcap X_s,W)
$$ is well defined. It is called the {\it equivalence} of $Z$ in 
$X_1\dcap\cdots\dcap X_s$ (cf. \cite[9.1]{fulton}). 

In analogy with \cite{twor}, one can define a local  variant of
this number as follows. For every $p$,
$$
\sum a_i\mult_pZ_i({\cal L})
$$ is constant on a Zariski open subset of the ${\cal L}$-s.
 I denote it by
$$
\mult_p (X_1\dcap\cdots\dcap X_s).
$$ There is an inequality
$$
\mult_p(X_1\dcap\cdots\dcap X_s)\leq 
\sum_{p\in W\in B(X_1\cap\cdots\cap X_s)}\deg (X_1\dcap\cdots\dcap
X_s,W).
$$
\end{defn}

We need to set up a correspondence between ideal sheaves and
algebraic cycles. This does not work as well as the usual
correspondence between subschemes and ideal sheaves, but it is
better suited for our purposes. Another way of going from cycles to
ideal sheaves is studied in  section 4. 

\begin{defn}\label{ass-id.def}
 Let $X$ be a scheme and
$Z=\sum a_i[Z_i]$ an effective cycle.  Let $I(Z_i)\subset \o_X$
denote the ideal sheaf of $Z_i$. Define the
 {\it ideal sheaf of $Z$} by
$$ I(Z):=\prod_iI(Z_i)^{a_i}\subset \o_X.
$$ It is clear that $I(Z_1+Z_2)=I(Z_1)I(Z_2)$.
\end{defn}

\begin{defn}\label{ass-cyc.def} Let $F$ be any coherent sheaf on
$X$ and $F_i\subset F$  the subsheaf of sections whose support has
codimension at most
$i$.  Let $x_{ij}$ be the generic points of the irreducible
components $X_{ij}\subset \supp(F_i/F_{i-1})$. Set
$$ Z(F):=\sum_{ij} (\len_{x_{ij}}F_i)\cdot [X_{ij}].
$$
$Z(F)$ is called the {\it cycle associated to $F$}.

Let $Q_{ij}\subset \o_X$ be the ideal sheaf of $X_{ij}$ and
$b_{ij}:=\len_{x_{ij}}F_i$. Then 
$\prod_jQ_{ij}^{b_{ij}}$ maps $F_i$ to $F_{i-1}$, thus
$I(Z(F))\subset \ann(F)$.

In particular, if $J\subset \o_X$ is an ideal sheaf  then
$$ I(Z(\o_X/J))\subset J.
$$

If $X$ is proper and $L$ is a line bundle on $X$ then one can define
the {\it $L$-arithmetic degree} of a sheaf $F$ by
$$
\tdeg_LF:=\deg_L Z(F). 
$$  If $I\subset \o_X$ is an ideal sheaf then the  arithmetic
degree of $\o_X/I$ is also called the arithmetic degree of $I$ and
denoted by $\tdeg_LI$. Note that there is a  possibility of
confusion since $I$ is also a sheaf.

This definition is  very natural and it appeared in several
different places (see, for instance,
\cite{harts, null, bay-mumf}).  The concept was  used extensively
in  many papers (cf. \cite{stv}). 
\end{defn}

\begin{lem}\label{I-prod.lem} Let $X_1,\dots,X_m$ be schemes and 
$Z_i$   a cycle on $X_i$ for every $i$. Let $\pi_i:\prod_j X_j\to
X_i$ be the $i$-th coordinate projection. Then
$$ I(\prod_j Z_j)\subset (\pi_1^*I(Z_1),\dots,\pi_m^*I(Z_m)).
$$
\end{lem}

Proof. Using induction, it is sufficient to prove the case $m=2$.
Let $Z_k=\sum_j a_{kj}Z_{kj}$, then
$$ Z_1\times Z_2=\sum_{ij}a_{1i}a_{2j}(Z_{1i}\times Z_{2j}).
$$ If $I,J$ are arbitrary ideals and $a,b\geq 1$, then 
$$ (I,J)^{ab}\subset (I,J)^{a+b-1}\subset (I^a,J^b).
$$ Using this on $X_1\times X_2$, we obtain that
\begin{eqnarray*} I(Z_1\times Z_2)&=& 
\prod_{ij}I(Z_{1i}\times Z_{2j})^{a_{1i}a_{2j}}\\ &=&
\prod_{ij}(\pi_1^*I(Z_{1i}), \pi_2^*I(Z_{2j}))^{a_{1i}a_{2j}}\\ 
&\subset &
\prod_{ij}(\pi_1^*I(Z_{1i})^{a_{1i}},
\pi_2^*I(Z_{2j})^{a_{2j}})\\  &\subset &
\prod_{j}(\pi_1^*\prod_{i}I(Z_{1i})^{a_{1i}},
\pi_2^*I(Z_{2j})^{a_{2j}})\\  &= &
\prod_{j}(\pi_1^*I(Z_1),
\pi_2^*I(Z_{2j})^{a_{2j}})\\  &\subset & (\pi_1^*I(Z_1),
\pi_2^*I(Z_2)).\qed
\end{eqnarray*}

\section{Integral closure of ideals}

In this section we recall some relevant facts concerning integral
closure of ideals. \cite[Chap.I]{teiss}  serves as a good 
 general reference.

\begin{defn}\label{intclos.def}
 Let $R$ be a ring and $I\subset R$ an ideal. 
 $r\in R$ is called {\it integral over $I$} if 
$r$ satisfies an equation
$$ r^k+\sum_{j=1}^k i_jr^{k-j}=0\qtq{where $i_j\in I^j$.}
$$ All elements integral over $I$ form an ideal $\overline{I}$,
called the {\it integral closure} of $I$. 

We use the following easy properties of the integral closure.
\begin{enumerate}
\item $(\overline{I}, \overline{J})\subset \overline{(I,J)}$,
\item $\overline{I_1}\cdot\overline{I_2}\subset \overline{I_1I_2}$, 
and so $(\overline{I})^m\subset \overline{I^m}$.
\end{enumerate}
\end{defn}

We also need the following special case of the  Brian\cedilla
con--Skoda theorem. A short proof of it  can be found in
\cite[p.101]{li-te}.

\begin{thm}\cite{bri-sko}\label{bri-sko} If $R=K[x_1,\dots,x_n]$
(or more generally, if $R$ is regular of dimension $n$) then
$\overline{I^n}\subset I$.\qed
\end{thm}

The following result gives the best way to compare integral
closures (cf. \cite[I.1.3.4]{teiss}). 

\begin{thm}[Valuative criterion of integral
dependence]\label{val-crit.thm}
 Let $R$ be a  ring and $I, J\subset R$ two ideal. The following
are equivalent.
\begin{enumerate}
\item  $J\subset \bar I$.
\item  If $p:R\to S$ is any homomorphism of $R$ to a DVR $S$  then
$p(J)\subset p(I)$.
\end{enumerate} If $K$ is an algebraically closed field and $R$ a
finitely generated
$K$-algebra then in (2) it is sufficient to use homomorphisms to
the power series ring 
$K[[t]]$. \qed
\end{thm}

Integral closures usually do not commute with taking quotients, but
this holds in some special cases.

\begin{lem}\label{intcl.rest.lem} Let $I\subset K[x_1,\dots,x_n]$
be an ideal. Then
$$
\overline{(I,x_n)}/(x_n)= \overline{(I,x_n)/(x_n)}.
$$
\end{lem}

Proof. If $J_1\subset J_2\subset R$ are ideals then
$\overline{J_2}/J_1\subset \overline{J_2/J_1}$ always holds using 
(\ref{intclos.def}).  If $R\to R/J_1$  splits (as a ring
homomorphism) then any equation over $R/J_1$ can be lifted to an
equation over $R$, showing the other containment.\qed
\medskip

We need two lemmas about ideals given by algebraic families of
generators.

\begin{lem} \label{zdense.gens.lem} Let $K$ be an infinite field,
$R$  a $K$-algebra and $L\subset R$ a finite dimensional
$K$-vectorspace. Let $U$ be a $K$-variety and 
$$ F:U\to L\qtq{given by} u\mapsto r_u
$$ a $K$-morphism. Let $V\subset U$ be Zariski dense. Then there is
an equality of ideals
$$ (r_u:u\in V)=(r_u:u\in U).
$$
\end{lem}

Proof. If $J$ is any ideal in $R$ then $L\cap J$ is  a sub vector
space in $L$. Thus $\{u\in U:r_u\in  J\}$ is  Zariski closed in $U$.
Set $J=(r_u:u\in V)$.  Since $V$ is dense in $U$, we obtain that
$r_u\in J$ for every $u\in U$.\qed

\begin{lem} \label{zdense.prod.lem} Notation as in
(\ref{zdense.gens.lem}).  Assume in addition that $U$ is
irreducible. Let
$u\mapsto r_u$ and $u\mapsto s_u$ be 
$K$-morphisms from $U$ to $L$. Let $V\subset U$ be Zariski dense.
Then we have an equality of ideals
$$
\overline{(r_us_u:u\in V)}=
\overline{(r_u:u\in U)\cdot (s_u:u\in U)}.
$$
\end{lem}

Proof. Let $p:R\to S$ be any homomorphism to a DVR. An ideal in $S$ is
characterized by the minimum order of vanishing of its elements. We
need to prove that both ideals above give the same number.

The order of vanishing of each $p(r_u)$ in $S$ is a lower semi
continuous function of
$U$, thus it achieves the minimum value on a  dense open subset of
$U$.  Similarly for $p(s_u)$. 
 Thus we can choose a  $u\in V$ where 
 both $p(r_u)$  and $p(s_u)$  achieve their minimum. \qed

\begin{exmp}   Let $R=K[x,y], L=\{ax+by\}, U=K, r_u=x-uy,
s_u=x+uy$.  Then $(r_u:u\in U)\cdot (s_u:u\in U)=(x,y)^2$ is
different from $(r_us_u:u\in U)=(x^2,y^2)$.  This  shows that 
(\ref{zdense.prod.lem}) fails without integral closure.

Another such  example is given in (\ref{axes.exmp}).
\end{exmp}

\begin{rem} It is easy to check that the ideals $(I_n,s)$ in
(\ref{algbez.false.exmp}) are integrally closed, hence
 integral closure alone cannot remove the embedded primes, even in a
geometrically very simple situation. 
\end{rem}

\section{The ideal of Chow  equations}

Let $K$ be a field and $Z$  any effective cycle in $\a^n$.  In this
section  we define an ideal in $K[x_1,\dots,x_n]$, called the ideal
of Chow  equations of
$Z$.   The main advantage of this notion is that it behaves well
with respect to arbitrary  hyperplane sections. This is the crucial
property that one needs for the applications. On the other hand, 
the ideal of Chow  equations is quite difficult to analyze and I
leave several basic questions unresolved. (The explanation of the
name and other variants are discussed in (\ref{why-chow?}).)

\begin{defn}\label{con-id.def} 
 Let $Z=\sum a_iZ_i$ be a  purely $d$-dimensional cycle in $\a^n$. 
Let $\pi:\a^n\to \a^{d+1}$ be a projection such that $\pi:Z_i\to
\a^{d+1}$ is finite for every $i$.  We call such a projection {\it
allowable}.

The center of the projection $\pi$ is a linear space
$L\subset \p^n\setminus \a^n$ of dimension $n-d-2$ and $\pi$ is
allowable iff $L$  is disjoint from
$\cup_i\bar Z_i$. This shows that allowable projections can be
parametrized by an irreducible  quasiprojective variety.

If $\pi$ is allowable then $\pi_*(Z)$ is a well defined codimension
1 cycle in
$\a^{d+1}$, and so it corresponds to a hypersurface. Choose an
equation of this hypersurface and pull it back by $\pi$ to obtain a
polynomial $f(\pi,Z)$.

Assume first that $K$ is infinite.  Define the {\it ideal of Chow 
equations} of $Z$  in the polynomial ring $K[\a^n]\cong
K[x_1,\dots,x_n]$  as
$$ I^{ch}(Z):=(f(\pi,Z): \mbox{ $\pi$ is allowable})\subset K[\a^n].
$$ For technical reasons we   frequently work with  the integral
closure of this ideal, denoted by $\overline{I^{ch}}(Z)$.

We see in (\ref{coneq.indep.cor}) that these are independent of the
base field. Thus if $K$ is finite, one can define 
$I^{ch}(Z)$ by taking any infinite field extension of $K$ first.

Finally, if $Z=\sum a_iZ_i$ is any effective cycle then  write $Z$
as a sum   $Z=\sum Z^d$ where $Z^d$ has pure dimension
$d$ and  set
$$ I^{ch}(Z)=\prod_d I^{ch}(Z^d).
$$ Its integral closure is denoted by $\overline{I^{ch}}(Z)$.  A
product formula in terms of the $Z_i$ is given in  (\ref{prod-lem}),
but this only works for the integral closures.
\end{defn}

\begin{rem}\label{why-chow?}

The ideals $I^{ch}(Z)$ were first considered by \cite{cayley} and 
$I^{ch}(Z)$ is  essentially equivalent to  the 
 Chow form of $Z$, as explained in
\cite{catan, ds95}. This equivalence clarifies the definition of
$I^{ch}(Z)$, but it obscures other   versions  of this concept.

In (\ref{con-id.def}) we consider {\it linear} projections
$\pi:\a^n\to \a^{d+1}$. It is, however, possible to use larger
classes of morphisms. For instance we can allow $\pi$ to be any
algebraic automorphism of $\a^n$ followed by a projection or we can
even allow $\pi$ to be any smooth morphism.  The latter case can be
 localized in various topologies.

More generally, if $R$ is any smooth $K$-algebra and $Z$  a
$d$-cycle on $\spec R$ then  one can define the  ideal of locally
Chow  equations (using \'etale or analytic topology or even working
formally) and these ideals behave well with respect to
intersections with smooth divisors. Here  I concentrate on the
simpler case of linear projections. I was unable to decide if the
various  definitions  give the same ideals for a cycle in $\a^n$.
\end{rem}

\begin{say} I do not know if it is essential to consider the
integral closure or not in the definition above.  The examples
(\ref{noncl.id.exmp}, \ref{charp.exmp}) show that
$I^{ch}(Z)$ is not integrally closed in general. More importantly,
the crucial property  (\ref{prod-lem}) fails without integral
closure as shown by
 (\ref{axes.exmp}). 

The main question is whether 
 (\ref{main.tech.thm}) holds without integral closure on the right
hand side.  This would eliminate the extra factor $(n+1)$ in
(\ref{eff.null.intr.thm}).  I do not know the answer.
This question is related to the degree bounds considered in
\cite[Sec.\ 4]{sturmf}.
\end{say}

As a special case of (\ref{zdense.gens.lem}) we obtain:

\begin{lem}\label{zdense-lem}
 Let $\{\pi_{\lambda}:\lambda \in \Lambda\}$ be a Zariski dense set
of allowable projections as in (\ref{con-id.def}).  Then
$$ I^{ch}(Z)=(f(\pi_{\lambda},Z)\vert \lambda \in
\Lambda).\qed
$$
\end{lem}

\begin{cor}\label{coneq.indep.cor}
 $I^{ch}(Z)$ is independent of the base field $K$. That is, if
$L\supset K$ is a field extension, then
$$ I^{ch}(Z)\otimes_KL=I^{ch}(Z_L).
$$
\end{cor}

Proof. If  $K$ is infinite, then the projections defined over $K$
form a  Zariski dense set of the  projections defined over $L$.
Thus by (\ref{zdense-lem}) we obtain the same ideals.

For finite $K$  we defined $I^{ch}(Z)$ by forcing the above formula
to hold.\qed

\begin{exmp} Let $X\subset \a^n$ be a smooth subvariety with ideal
sheaf $I(X)$.  Then $I^{ch}(X)=I(X)$ and $I^{ch}(a\cdot X)=I(X)^a$.
More generally, let  $Z=\sum a_iZ_i$ be any cycle. Then the above
relationship holds near any smooth point of $\supp Z$. 

Thus $I^{ch}(Z)$ is interesting only near the singular points of
$\supp Z$.

Let $p\in \supp Z$ be a point of multiplicity $d$ and $m_p$ the
ideal of
$p$. A general projection $\pi(Z)$ has multiplicity $d$ at $\pi(p)$,
thus each $f(\pi,Z)$ has  multiplicity $\geq d$ at $p$.  This shows
that $I^{ch}(Z)\subset m_p^d\cap I(Z)$. 

By \cite[1.14]{catan}, if $Z$ has codimension at least 2 then
$I(Z)\neq I^{ch}(Z)$ along the singular locus of $Z$. 
\end{exmp}

\begin{exmp}\label{low.embdim.exmp} Let $\a^{n-k}\subset \a^n$ be
the   subspace
$(x_n=\cdots=x_{n-k+1}=0)$. Let $Z$ be a $d$-cycle on $\a^{n-k}$
and $j_*Z$ the corresponding cycle on $\a^n$. We would like to
compare $I^{ch}(Z)$ and $I^{ch}(j_*Z)$. 

A general projection of $j_*Z$ can be obtained as a projection
$\rho:\a^n\to \a^{n-k}$ followed by a general projection
$\pi:\a^{n-k}\to \a^{d+1}$. This shows that
$$ f(\pi\circ \rho, j_*Z)=f(\pi, Z)(x_1+L_1,\dots,x_{n-k}+L_{n-k}),
$$ where the $L_i$ are linear forms in $x_{n-k+1},\dots, x_n$
defining
$\rho$. 

This shows that the restriction map
$$ I^{ch}(j_*Z)\onto I^{ch}(Z)\qtq{is surjective.}
$$
\end{exmp}

\begin{exmp}\label{noncl.id.exmp}
 Let $X\subset \a^n$ be defined by equations
$g(x_1,\dots,x_{n-1})=x_n=0$. A general  projection of $X$ is 
isomorphic to $X$ and, at least in characteristic zero, 
\begin{eqnarray*}
f(\pi,X)&=&g(x_1+a_1x_n,\dots,x_{n-1}+a_{n-1}x_n)\\ &=& \sum_I
c_Ia^Ix_n^{|I|}\frac{\partial^Ig}{\partial x^I},
\end{eqnarray*} where the $c_I$ are nonzero constants and $a_i\in
K$. Since the $a_i$ can vary independently, we see that the
$f(\pi,X)$ generate the ideal
$$
\left(x_n^{|I|}\frac{\partial^Ig}{\partial x^I}:
I=(i_1,\dots,i_{n-1})\right).
$$

Consider for instance  the case $n=3$ and $g=x_1^3+x_2^5$. Then 
$$ I^{ch}(X)=(x_1^3+x_2^5,x_1^2x_3,x_1x_3^2,x_3^3,
x_2^4x_3,x_2^3x_3^2,x_2^2x_3^3,x_2x_3^4,x_3^5).
$$
 $x_2^2x_3^2$ is integral over $I^{ch}(X)$  (since
$(x_2^2x_3^2)^2-x_3^3\cdot x_2^4x_3=0$) but it is not in
$I^{ch}(X)$. Hence $\overline{I^{ch}}(X)\neq I^{ch}(X)$. 
\end{exmp}

\begin{exmp}\label{charp.exmp}
 Assume that $K$ has characteristic $p$ and let $0\in \a^2$ be the
origin with ideal $(x,y)$. Set $Z=p[0]$.  Then $f(\pi,Z)=(ax+by)^p$
for some $a,b$, thus
 $I^{ch}(Z)=(x^p,y^p)$. Its integral closure is  the much bigger
ideal $(x,y)^p$.
\end{exmp}

The next lemma gives a product formula for $\overline{I^{ch}}(Z)$.
This result is crucial for the applications and it fails if we do
not take the integral closure, as the example after the lemma shows.

\begin{lem}\label{prod-lem} Let $Z=\sum a_iZ_i$ be an effective 
cycle. Then
$$
\overline{I^{ch}}(Z)=\overline{\prod_i I^{ch}(Z_i)^{a_i}}.
$$
\end{lem}

Proof.  It is enough to check this for pure dimensional cycles.

Let $\pi$ be any allowable projection for $Z$. Then $\pi$ is
allowable for every $Z_i$ and
$f(\pi,Z)=\prod_if(\pi,Z_i)^{a_i}$ which proves the containment
$\subset$. The converse follows by a repeated application of
(\ref{zdense.prod.lem}).\qed

\begin{exmp}\label{axes.exmp}
 Choose $n\geq 3$ odd and in $\a^n$ consider the 1--cycle
 of the $n$  coordinate axes $Z=\sum_{i=1}^n Z_i$. Then
$$
\prod_i I^{ch}(Z_i)=\prod_i(x_1,\dots,\widehat{x_i},\dots,x_n).
$$ As a vectorspace this has a basis consisting of all monomials of
degree at least $n$ which involve at least 2 variables.

On the other hand, I claim that $I^{ch}(Z)$ does not contain the
monomial  $x_1\cdots x_n$. 

 Let $\pi:(x_1,\dots,x_n)\mapsto (\sum a_ix_i,\sum b_ix_i)$ be a
projection.     This gives the equation
$$ f(\pi,Z)=\prod_i\left(\sum_j m_{ji}x_j\right)
\qtq{where $m_{ji}=a_jb_i-a_ib_j$.}
$$ Thus twice the coefficient of the  $x_1\cdots x_n$ term is

\begin{eqnarray*} 2\cdot \sum_{\sigma\in S_n}\prod_j m_{j\sigma(j)}
&=&\sum_{\sigma\in S_n}\left(\prod_j m_{j\sigma(j)}+\prod_j
m_{j\sigma^{-1}(j)}\right)\\ &=&\sum_{\sigma\in S_n}\left(\prod_j
m_{j\sigma(j)}+\prod_j m_{\sigma(j)j}\right)\\ &=&\sum_{\sigma\in
S_n}\left(\prod_j m_{j\sigma(j)}+(-1)^n\prod_j
m_{j\sigma(j)}\right)=0.
\end{eqnarray*}

If $n=3$ then it is easy to compute that 
$\prod_i I^{ch}(Z_i)=(I^{ch}(Z),x_1x_2x_3)$. I have not checked what
happens for $n\geq 5$ or for even values of $n$. 
\end{exmp}

The following   result of \cite[Thm.B]{amor} shows that
$\overline{I^{ch}}(Z)$ contains a fairly small power of $I(Z)$.
(The statement in \cite{amor} is slightly different since he is
working with $I^{ch}(Z)$, but his proof actually gives this
version.) 

\begin{thm}\label{jks.lem}\cite{amor}
 Let $Z=\sum a_iZ_i$ be a cycle in $\a^n$. Let $I(Z_i)$ denote the
ideal of $Z_i$. Then
$$
\prod_iI(Z_i)^{a_i\deg Z_i}\subset \overline{I^{ch}}(Z).
$$ More precisely, if $x\in \a^n$ is a point then
$$
\prod_iI(Z_i)^{a_i\mult_xZ_i}\subset \overline{I^{ch}}(Z)
\qtq{in a neighborhood of $x$.}\qed
$$
\end{thm}

\begin{lem}\label{con<ord.lem} Let $Z$ be a cycle in $\a^n$. Then
$$ I^{ch}(Z)\subset I(Z).
$$ Let $J\subset K[x_1,\dots,x_n]$ be an ideal. Then
$$ I^{ch}(Z(K[x_1,\dots,x_n]/J))\subset J.
$$
\end{lem}

Proof. Because of the multiplicative definitions of $I(Z)$
(\ref{ass-id.def}) and of $I^{ch}(Z)$ (\ref{con-id.def}), it is
sufficient to prove the first claim in case $Z$ is pure dimensional.

Write $Z=\sum a_iZ_i$ and let $\pi:\a^n\to \a^{d+1}$ be an allowable
projection. Then $\pi_*[Z]=\sum a_i\pi_*[Z_i]$, so
$f(\pi,Z)=\prod f(\pi,Z_i)^{a_i}$. $f(\pi,Z_i)\in I(Z_i)$, so
$f(\pi,Z)\in \prod I(Z_i)^{a_i}=I(Z)$.

The second part follows from the first and from
(\ref{ass-cyc.def}).\qed

\section{The ideal of Chow  equations and intersection theory}

The next result is the key property of the ideal of Chow  equations.

\begin{lem}\label{coneq.main.lem}
 Let $X\subset \a^n$ be an irreducible and reduced subvariety and
$H=(x_n=0)$  a hyperplane not containing $X$. Then
$$ I^{ch}(X\ncap H)\subset (I^{ch}(X),x_n).
$$
\end{lem}

Proof.  By (\ref{coneq.indep.cor}) we may assume that the base
field is infinite. Choose a general linear subspace $L\subset
\bar H \setminus H \setminus\bar X$ of dimension $n-d-2$.
$\dim (\bar H \setminus H)=n-2$ and $\dim (\bar X\cap \bar H)\leq
d-1$. Since $(n-d-2)+(d-1)<n-2$, $L$ is disjoint from $\bar X$. 
Let $\pi':H\to H'$ and $\pi:\a^n\to \a^{d+1}$ be the  projections
with center $L$. Let $\rho:\a^n\to H$ be a projection and set
$\pi'':=\pi\circ \rho: \a^n\to H'$. The 3 projections  appear in
the following diagram:
$$
\begin{array}{cll} H & \subset & \a^n\\
\pi'\downarrow{\ \ } & \swarrow \pi'' & \downarrow\pi\\ H'& \subset
& \a^{d+1}
\end{array}
$$
$ X\ncap H$ can be viewed as a cycle on $H$; in such a case I
denote it by $Z$. 

$\pi_*(X)$ is a hypersurface in $\a^{d+1}$ and $\pi'_*(Z)$ is a
hypersurface in $H'$ such that $\pi'_*(Z)=\pi_*(X)\cap H'$.   Thus 
$$ f(\pi',Z)(x_1,\dots,x_{n-1})=f(\pi,X)(x_1,\dots,x_{n-1},0).
$$

As in (\ref{low.embdim.exmp}),   the generators of $I^{ch}(X\ncap
H)$ are of the form
\begin{eqnarray*} f(\pi'', X\ncap H)&=& 
f(\pi',Z)(x_1+a_1x_n,\dots,x_{n-1}+a_{n-1}x_n)\\ &\equiv
&f(\pi',Z)(x_1,\dots,x_{n-1})\mod (x_n)\\ &\equiv
&f(\pi,X)(x_1,\dots,x_{n-1}, x_n)\mod (x_n). \qed
\end{eqnarray*}
\medskip

\begin{rem}\label{coneq.main.lem-rem}
 More generally, (\ref{coneq.main.lem}) also holds if $X$ is a pure
dimensional cycle and $H$ does not contain any of its irreducible
components.
\end{rem}

The generalization to intersecting with several linear equations is
formal, but the induction seems to require the use of integral
closure, as shown by the following example. The final result
itself, however, may not need integral closure.

\begin{exmp} The 3 coordinate axes 
$C_x,C_y,C_z$ in $\a^3$ can be defined by the determinantal
equations
$$
\rank
\left(
\begin{array}{ccc} x&y&0\\ z&y&z
\end{array}
\right)
\leq 1.
$$
$C_x+C_y+C_z$ is a hyperplane section of the surface 
$$ Z_1 \qtq{given by equations}
\rank
\left(
\begin{array}{ccc} x&y+as&bs\\ z&y+cs&z+ds
\end{array}
\right)
\leq 1.
$$ For general $a,b,c,d$, $Z_1$ is the cone over a rational normal
curve. By explicit computation,
$xyz\in I^{ch}(Z_1)$. As we remarked in (\ref{axes.exmp}), this
implies that
$$ I^{ch}(C_x)\cdot I^{ch}(C_y)\cdot I^{ch}(C_z)
\subset (I^{ch}(Z_1),s).
$$ Next consider the surface
$$ Z_2 \qtq{with equations}
\rank
\left(
\begin{array}{ccc} x&y+as^2&bs^2\\ z&y+cs^2&z+ds^2
\end{array}
\right)
\leq 1.
$$ For general $a,b,c,d$ this defines a rational triple point.  By
explicit computation,
$xyz\not\in (I^{ch}(Z_2),s)$, which implies that
$$ I^{ch}(C_x)\cdot I^{ch}(C_y)\cdot I^{ch}(C_z)
\not\subset (I^{ch}(Z_2),s).
$$
\end{exmp}

\begin{lem} \label{surfcut.ind} Let 
$Z$ be a cycle on $\a^n$ and $H_i=(\ell_i=0)$ hyperplanes.
 Then
 $$ I^{ch}(Z\ncap H_1\ncap \cdots\ncap H_m)\subset
\overline{(I^{ch}(Z),\ell_1,\dots,\ell_m)}.
$$
\end{lem}

Proof. 
 Consider first the case when $Z$ is irreducible and $m=1$. The
claim is trivial if $Z\subset H_1$ and the $Z\not\subset H_1$ case
is treated in (\ref{coneq.main.lem}).

Next we prove the $m=1$ case by  induction on the number of
irreducible components of $Z$.
\begin{eqnarray*} I^{ch}((Z_1+Z_2)\ncap \ell_1)&\subset & 
\overline{I^{ch}(Z_1\ncap \ell_1)\cdot I^{ch}(Z_2\ncap \ell_1)}
                   \qtq{(by (\ref{prod-lem}))}\\ &\subset &
\overline{(I^{ch}(Z_1),\ell_1)\cdot (I^{ch}(Z_2),\ell_1)}
\qtq{(by induction)}\\ &\subset &
\overline{(I^{ch}(Z_1)I^{ch}(Z_2),\ell_1)}\\
 &\subset &
\overline{(I^{ch}(Z_1+Z_2),\ell_1)}\qtq{(by (\ref{prod-lem}))}
\end{eqnarray*}

Finally the  case $m>1$  is established by induction using the
chain of inclusions
\begin{eqnarray*} I^{ch}(Z\ncap \ell_1\ncap  \ell_2)&\subset &
\overline{(I^{ch}(Z\ncap\ell_1),\ell_2)}\\ &\subset&
\overline{(\overline{(I^{ch}(Z), \ell_1)},\ell_2)}\\ &\subset&
\overline{(I^{ch}(Z),\ell_1,\ell_2)}.\qed
\end{eqnarray*}

We are ready to formulate our main technical theorem.

\begin{thm}\label{main.tech.thm}  Let $K$ be a field and
$Z_1,\dots,Z_m$    cycles in $\a^n$.  Let
$Z_1\ncap\cdots\ncap Z_m$ be any of the
 intersection cycles. Then
$$ I^{ch}(Z_1\ncap\cdots\ncap Z_m) \subset 
\overline{(I(Z_1),\dots,I(Z_m))}.
$$
\end{thm}

Proof.  Choose an identification
$\a^{nm}=\a^n\times\cdots\times\a^n$ ($m$-times) and let
$\Pi_r:\a^{nm}\to \a^n$ be the projection onto the $r$-th factor.
Let $\Delta\subset \a^n\times\cdots\times\a^n$ denote the
 diagonal.

Choose   an ordered set of hyperplanes  ${\cal L}:=(L_i=(\ell_i=0):
i=1,\dots,n(m-1))$   in
$\a^{nm}$  whose intersection is $\Delta$. This gives us a cycle
$(Z_1\ncap\cdots\ncap Z_m,{\cal L})$ which we view as a cycle in
$\a^{nm}$. 

Applying (\ref{surfcut.ind}) we obtain that
$$ I^{ch}(Z_1\ncap\cdots\ncap Z_m,{\cal L})
\subset 
\overline{(I^{ch}(Z_1\times\cdots\times Z_m),
\ell_1,\dots,\ell_{n(m-1)})}.
$$
$I^{ch}(Z_1\times\cdots\times Z_m)
\subset I(Z_1\times\cdots\times Z_m)$ by (\ref{con<ord.lem}) and
using
 (\ref{I-prod.lem}) this gives the inclusion
$$ I^{ch}(Z_1\ncap\cdots\ncap Z_m,{\cal L})
\subset 
\overline{(\Pi_1^*I(Z_1),\dots, \Pi_m^*I(Z_m),
\ell_1,\dots,\ell_{n(m-1)})}.
$$ Let us restrict  to $\Delta$.  The left hand side becomes 
$I^{ch}(Z_1\ncap\cdots\ncap Z_m)$ by (\ref{low.embdim.exmp}), and
the right hand side becomes 
$\overline{(I(Z_1),\dots,I(Z_m))}$ by (\ref{intcl.rest.lem}). 
\qed

\begin{rem} It is possible that (\ref{main.tech.thm}) can be
considerably sharpened. The strongest and most natural statement
would be
$$ I^{ch}(Z_1\ncap\cdots\ncap Z_m) \subset 
(I^{ch}(Z_1),\dots,I^{ch}(Z_m)).
$$ For the applications the main point would be to get rid of the
integral closure since this would eliminate the extra factor 
$(n+1)$ in (\ref{eff.null.intr.thm}). 
\end{rem}

\section{Effective Nullstellensatz}

We are ready to formulate and prove the precise technical versions
of our main theorems, using the notion of arithmetic degree as defined
in (\ref{ass-cyc.def}).

\begin{thm}[Algebraic B\'ezout theorem]\label{bez.thm}
  Let $K$ be any field and \linebreak 
$I_1,\dots,I_m$
   ideals in $K[x_1,\dots,x_n]$.  Then there are   prime ideals
$P_j\supset (I_1,\dots,I_m)$ and natural numbers
$a_j$ such that
\begin{enumerate}
\item  $\prod_j P_j^{a_j}\subset (I_1,\dots,I_m)$, and
\item $\sum_j a_j\leq n\cdot \prod_i\tdeg I_i$.
\end{enumerate}
\end{thm}

Proof. As in (\ref{ass-cyc.def}), set $Z_i=Z(I_i)$ and let 
$\sum b_jX_i=(Z_1\ncap\cdots\ncap Z_m,{\cal L})$ be any of the
intersection cycles defined in (\ref{vt-cyc.def}). Set $d_j:=\deg
X_j$, then 
$\sum_j b_jd_j\leq \prod_i \deg Z_i$ by (\ref{bezout.thm}).  By
(\ref{main.tech.thm}),
$$
\prod_jI^{ch}(X_j)^{b_j} \subset 
\overline{(I(Z_1),\dots,I(Z_m))}.
$$
$I(X_j)^{d_j}\subset \overline{I^{ch}}(X_j)$ by (\ref{jks.lem}),
and so we obtain that
$$
\prod_jI(X_j)^{b_jd_j} \subset 
\overline{(I(Z_1),\dots,I(Z_m))}.
$$
$I(Z_s)\subset I_s$ by (\ref{con<ord.lem}), hence
$$
\overline{(I(Z_1),\dots,I(Z_m))}\subset \overline{(I_1,\dots,I_m)}.
$$
$(I_1,\dots,I_m)^n\subset \overline{(I_1,\dots,I_m)}$ by
(\ref{bri-sko}). Putting these together we get that
$$
\prod_jI(X_j)^{nb_jd_j} \subset  (I_1,\dots,I_m).
$$
 Setting $P_j:=I(X_j)$  and $a_j:=nb_jd_j$ gives
(\ref{bez.thm}).\qed

\begin{thm}[Effective Nullstellensatz]\label{eff.null.thm} Let $K$
be any field and \linebreak
$I_1,\dots,I_m$    ideals in $K[x_1,\dots,x_n]$.  The following are
equivalent:
\begin{enumerate}
\item  $I_1,\dots,I_m$ have no common zero in $\bar K^n$.
\item There are polynomials $f_j\in I_j$  such that
$$
\sum_j f_j=1\qtq{and} \deg f_j
\leq (n+1)\cdot 
\prod_i\tdeg I_i.
$$
\end{enumerate}
\end{thm}

Proof.  It is clear that $(2)\Rightarrow (1)$. To see the
converse,  introduce a new variable $x_0$ and let
$\tilde I_s\subset K[x_0,\dots,x_n]$ denote the homogenization of
$I_s \subset K[x_1,\dots,x_n]$. Then $\tdeg \tilde I_s=\tdeg I_s$
and $x_0$ is contained  in the radical of
$(\tilde I_1,\dots, \tilde I_m)$, hence it is contained in any
prime ideal containing  $(\tilde I_1,\dots, \tilde I_m)$.  By
(\ref{bez.thm}) there are prime ideals $P_j$ and natural numbers
$a_j$ such that 
\begin{enumerate}
\item  $\prod_j P_j^{a_j}\subset (\tilde I_1,\dots,\tilde I_m)$, and
\item $\sum_j a_j\leq (n+1)\cdot \prod_i\tdeg I_i$.
\end{enumerate} Since $x_0\in P_j$ for every $j$, we see that
$$ x_0^{\sum a_j}\in (\tilde I_1,\dots,\tilde I_m).
$$ Thus there are $f_i\in I_i$ with homogenizations
$\tilde f_i$  such that
$$ x_0^{\sum a_j}=\sum_i  \tilde f_i
\qtq{and} \deg \tilde f_i= \sum_j a_j.
$$ Setting $x_0=1$ we obtain  (\ref{eff.null.thm}).\qed

\section{\L ojasiewicz inequalities}

Next we turn to  applications of  these results to the study of 
\L ojasiewicz inequalities  and  separation exponents. These
results are essentially  reformulations  of \cite{cygan, ckt}. 

\begin{defn}\label{loj.ineq} Let $f$ be a real analytic function on
$\r^n$ and $Z:=(f=0)$. Fix a norm on $\r^n$ and set
$\dist(Z,x):=\inf_{z\in Z}\norm{x-z}$. 
\cite[p.124]{loj-stud} proved that  for every compact set $K$ there
are
$m, C>0$ such that
$$
\dist(Z,x)^m\leq C\cdot |f(x)|\qtq{for $x\in K$.}
$$ Any inequality of this type is  called a  {\it \L ojasiewicz
inequality}. 
\end{defn}

In general it is rather difficult to obtain an upper bound  for $m$
in terms of  other invariants of $f$. The problem becomes easier if
$\r^n$ is replaced by $\c^n$, but even in this case it is not
straightforward to obtain sharp upper bounds for $m$.  The question
was investigated in
\cite{brownaw-loc} and \cite{jks}.  Instead of $\c$, one can work
over any  algebraically closed field with an absolute value.

\begin{notation}
 Let $K$ be a field with an absolute value $|\ |$.  (The case when
$K=\c$ and $|\ |$ is the usual absolute value is the most
interesting, but the cases when 
$K$ is of positive characteristic or 
$|\ |$ is nonarchimedian are also of interest.)
 $|\ |$ induces a norm on $K^n$ by $\norm{{\bold
x}}:=(|x_1|+\cdots+|x_n|)^{1/2}$. This defines a distance on $K^n$
as in (\ref{loj.ineq}). 
\end{notation}

\begin{defn}\label{gen.int.def} Let $X$ be any topological space
and $F,G$ two sets of $K$-valued functions on $X$. We say that  $F$
is {\it integral over $G$}, denoted by
$F\ll G$, if the following condition holds:
\begin{enumerate}
\item[$(*)$] For every $f\in F$ and $x\in X$ there are
$g_1,\dots,g_m\in G$ and a constant $C$ such that $|f(x')|\leq
C\max_i|g_i(x')|$ for every
$x'$ in a neighborhood of $x$. 
\end{enumerate}

If $F$ and $G$ are continuous (which will always be the case for us)
then $(*)$ is automatic if $g(x)\neq 0$ for some $g\in G$. Thus
$(*)$ is a local growth condition near the common zeros of $G$.
\end{defn}

The two notions of integral dependence are closely related:

\begin{lem}\label{int.same.lem} (cf. \cite[1.3.1]{teiss})
 Let
 $K$ be an algebraically closed field with an absolute value $|\ |$.
Let
$X$ be an affine  variety over $K$  (with the metric topology) and
$I\subset
\o_X$  an ideal sheaf. A polynomial function is  integral over $I$ 
in the sense of (\ref{intclos.def}) iff it is 
 integral over $I$  in the sense of (\ref{gen.int.def}).\qed 
\end{lem}

The relationship between the distance function and the ideal of
Chow  equations was established in earlier papers.

\begin{lem} \label{dist-id.lem} (cf. \cite[8]{jks},
\cite[3.7]{cygan})
 Let
 $K$ be an algebraically closed field  with an absolute value $|\
|$. Let $Z\subset K^n$ be an irreducible subvariety, $z\in Z$ a
point and $m=\mult_zZ$. Then, in a neighborhood of $z$, 
$$
\dist(Z,x)^m\ll I^{ch}_Z\ll I_Z\ll \dist(Z,x).\qed
$$
\end{lem}

The  main result about \L ojasiewicz inequalities is the following.

\begin{thm}[\L ojasiewicz inequality]\label{lojas.thm} (cf.\
\cite{ckt})  Let
 $K$ be an algebraically closed field with an absolute value $|\ |$.
Let $I_1,\dots,I_m $ be ideals in $K[x_1,\dots,x_n]$ and
 $X_1,\dots,X_m\subset K^n$ the corresponding  subschemes. Set
$D:=\tdeg I_1\cdots\tdeg I_m$.  Then 
$$
\dist(X_1\cap \cdots \cap X_m,z)^D\ll (I_1,\dots,I_m)\ll
\max_i\{\dist(X_i,z)\}.
$$
\end{thm}

Proof.  Let $(X_1\ncap\dots\ncap X_m)=\sum a_iZ_i$ be one of the
intersection cycles.
$\sum a_i\deg Z_i\leq D$ by (\ref{bezout.thm}) and $Z_i\subset
X_1\cap \cdots \cap X_m$ by construction. Thus
\begin{eqnarray*}
\dist(X_1\cap \cdots \cap X_m,z)^D&\ll&
\dist(X_1\cap \cdots \cap X_m,z)^{\sum a_i\deg Z_i}\\ &\leq &
\prod_i \dist(Z_i,z)^{a_i\deg Z_i}\\ & \ll &
\prod_i I^{ch}(Z_i)^{a_i}\qtq{(by (\ref{dist-id.lem}))}\\ & \subset
&\overline{(I_1,\dots,I_m)}\qtq{(by (\ref{bez.thm}))}\\  & \ll &
(I_1,\dots,I_m) \qtq{(by (\ref{int.same.lem}))}\\ &\ll & 
\max_i\{\dist(X_i,z)\}\qtq{(by (\ref{dist-id.lem})).\qed}
\end{eqnarray*}

With a similar proof we obtain the following
 local version.

\begin{cor}\label{local.sep.cor}(cf.\ \cite[4.5]{cygan}) 
 Let
 $K$ be an algebraically closed field with an absolute value $|\ |$.
Let $I_1,\dots,I_m $ be ideals in $K[x_1,\dots,x_n]$ and set
$D_p:=\mult_p(Z(I_1)\dcap\cdots\dcap Z(I_m))$.  Then, in a
neighborhood of
$p$, 
$$
\dist(X_1\cap \cdots \cap X_m,z)^{D_p}\ll (I_1,\dots,I_m)\ll
\max_i\{\dist(X_i,z)\}. \qed
$$
\end{cor}

\section{Application to deformation theory}

In usual deformation theory we are given a scheme $X_0$ and we 
would like to understand all flat families $\{X_t:t\in \Delta\}$
where
$\Delta$ is the unit disc. There are, however, some deformation
problems where we are interested in flat families $\{Y_t:t\in
\Delta\}$ where
$X_0=Y_0/(\mbox{embedded points})$, or, more generally, when $X_0$
and
$Y_0$  have the same fundamental cycles. This question arises for
instance in studying the Chow varieties. (See \cite{hp} or
\cite[Chap.I]{koll96} for definitions and properties of the Chow
varieties.) A point in the Chow variety of $\p^n$ is not a
subscheme but a pure dimensional cycle $W\in Z_d(\p^n)$. Thus if we
want to study the Chow variety near $W$ then we need to understand
the deformations of all subschemes $X\subset \p^n$ whose
fundamental cycle is
$W$. If $d\geq 1$ then there are infinitely many such subschemes
$X$ since adding embedded points does not change the fundamental
cycle.

Assume that we find a subscheme $X_0\subset \p^n$  whose
fundamental cycle is $W$ and a deformation $\{X_t:t\in \Delta\}$.
From the point of view of the Chow variety we are interested only
in the fundamental cycle of $X_t$ and not in $X_t$ itself.  Hence
the only case we need to study is when $X_t$ has no embedded points
for general $t$. The affine version of this problem can be stated as
follows.

\begin{question} Let $W\in Z_d(\a^n)$ be a $d$-cycle. Let $Y\subset
\a^n\times \a$ be a subscheme of pure dimension $(d+1)$ without
embedded points such that the second projection $\pi:Y\to \a^1$ is
flat.  Let $Y_0=\pi^{-1}(0)$ be the central fiber and  assume  that
$Z(Y_0)=W$. 

What can we say about $Y_0$ in terms of $W$?
\end{question}

This question is related to the problems considered in
\cite{flat-crit}.

As an application of (\ref{coneq.main.lem-rem}) we obtain the
following partial answer. This is a place where it would be  more
natural to use the ideal of locally Chow  equations (\ref{why-chow?}).

\begin{prop}\label{def.prop}
 With the above notation, $I^{ch}(W)\subset I(Y_0)$.\qed
\end{prop}

\begin{exmp} Consider the case when $W=[z=x^2-y^n=0]\in Z_1(\a^3)$.
As in (\ref{noncl.id.exmp}) we obtain that
$$ I^{ch}(W)=(x^2-y^n,z^2,xz, y^{n-1}z).
$$ On the other hand, in (\ref{algbez.false.exmp}) we found an
example of a deformation $S$ such that
$$ I(S_0)=(x^2-y^n,z^2,xz, y^{\frac{n-1}{2}}z).
$$ Using \cite{teiss-sn} we obtain that the length of
$I(W)/I(Y_0)$ is at most the arithmetic genus of $W$ which is
$\frac{n-1}{2}$. Comparing these two results we conclude that
$$ (x^2-y^n,z)\supset I(Y_0)\supset (x^2-y^n,z^2,xz,
y^{\frac{n-1}{2}}z)
$$ for every deformation $Y$. It is not hard to see that for every
ideal satisfying the above condition there is a corresponding
deformation.
\end{exmp}

Using (\ref{noncl.id.exmp}) one can compute $I^{ch}$ for all 
monomial plane curves in $\a^3$.  The results give strong
restrictions on $Y_0$ but I do not see how to get a complete answer
as in the above example. P.\ Roberts computed several examples of
monomial space curves and in each case $I^{ch}$ turned out to be
quite close to the ideal of the curve.

In higher dimensions $I^{ch}$ gives a very unsatisfactory answer
when $W=[X_0]$ where $X_0$ is normal.  By
\cite[III.9.12]{harts-book} in this case $Y_0=X_0$ but
$I^{ch}(X_0)\neq I(X_0)$ if $X_0$ is singular.

On the other hand, (\ref{def.prop})  gives information about $Y_0$
even if
$W$ has multiple components. I do not see how to get any 
information about $Y_0$ by other methods in this case.

\vskip1cm

\noindent University of Utah, Salt Lake City UT 84112 

\begin{verbatim}kollar@math.utah.edu\end{verbatim}

\end{document}